\documentclass[journal]{IEEEtran}

\ifCLASSINFOpdf
   \usepackage[pdftex]{graphicx}
\else
   \usepackage[dvips]{graphicx}
\fi

\usepackage[cmex10]{amsmath}
\usepackage{bm}

\usepackage{algorithmic}
\usepackage{algorithm}
\usepackage{array}
\usepackage{cancel}
\usepackage{cleveref}
\usepackage[dvipsnames]{xcolor}
\usepackage{setspace}

\usepackage{caption} %
\usepackage{subcaption} %

\allowdisplaybreaks

\def\removepowerplot{1}

\newfloat{model}{thp}{lop}
\floatname{model}{Model}

\newcommand{\edit}[1]{{\color{black}#1}}

\begin{document}
\bstctlcite{IEEEexample:BSTcontrol}
\title{Long Solution Times or Low Solution Quality: On Trade-Offs in Choosing a Power Flow Formulation for the Optimal Power Shutoff Problem
}

\author{\IEEEauthorblockN{Eric Haag$^*$,
Noah Rhodes$^*$, and 
Line Roald}
\IEEEauthorblockA{University of Wisconsin-Madison, Madison, Wisconsin, USA} %
}

\maketitle
\begin{abstract}
    The Optimal Power Shutoff (OPS) problem is an optimization problem that makes power line de-energization decisions in order to reduce the risk of igniting a wildfire, while minimizing the load shed of customers.  This problem, with DC linear power flow equations, has been used in many studies in recent years.  However, using linear approximations for power flow when making decisions on the network topology is known to cause challenges with AC feasibility of the resulting network, as studied in the related contexts of optimal transmission switching or grid restoration planning. This paper explores the accuracy of the DC OPS formulation and the ability to recover an AC-feasible power flow solution after de-energization decisions are made.  We also extend the OPS problem to include variants with the AC, Second-Order-Cone, and Network-Flow power flow equations, and compare them to the DC approximation with respect to solution quality and time. The results highlight that the DC approximation overestimates the amount of load that can be served, leading to poor de-energization decisions. %
    The AC and SOC-based formulations are better, but prohibitively slow to solve for even modestly sized networks thus demonstrating the need for new solution methods with better trade-offs between computational time and solution quality.
\end{abstract}
\begingroup
\def\thefootnote{}
\footnote{
\kern -3pt
\hrule width 0.4\columnwidth height 0.2pt
\kern 5pt
This work is supported by the U.S. National Science Foundation under the NSF CAREER award 2045860 and NSF ASCENT award 2132904.
}
\setcounter{footnote}{0}%
\endgroup
\def\thefootnote{*}\footnotetext{These authors contributed equally to this work}

\section{Introduction}
Climate change is expected to increase the risk of wildfires in many parts of the world~\cite{portner2022ipcc}. In recent years, several large and devastating fires ignited by power lines have increased scrutiny on utility practices to avoid wildfire ignitions. In California, utilities are routinely leveraging de-energization of power lines -- commonly referred to as public safety power shutoffs (PSPS) -- as a measure to prevent fires (and limit their liability) during periods of extreme risk. %
While PSPS is effective at reducing wildfire ignitions~\cite{pgeOct2019report}, it carries a significant cost in customer outages.
Prior work has considered how to optimally implement PSPS, a problem called \emph{Optimal Power Shutoff} (OPS), to balance the wildfire risk reduction with the resulting load shed~\cite{rhodes2020balancing}, reduce PSPS-caused power outages through energy storage~\cite{kody2022optimizing, astudillo2022managing} or microgrids~\cite{yang2022resilient}\cite{hanna2021optimal}, or transmission upgrades~\cite{taylor2022framework, bertoletti2022transmission, bayani2023resilient}.

A common aspect of prior work is that they all rely on a DC power flow approximation to reduce computational time. %
However, DC power flow does not account for reactive power flows and voltage constraints, and produces solutions that are not AC power flow feasible~\cite{baker2021solutions}. This is particularly true when the system is operating under stress, %
as can be expected when a substantial number of lines are switched off during a power shutoff.
Some data-driven methods for power management during periods of high wildfire risk use AC power flow, but only considers up to three de-energizations in the network~\cite{hong2022data}.

The goal of this paper is to assess how the choice of a power flow representation in the OPS problem impacts solution quality, solution time, and the trade-off between the two. While this is a question that is important for the OPS problem, it also arises in many other power system optimization problems which involve binary decisions regarding whether or not certain transmission lines should be in operation. 
Examples include transmission switching~\cite{hedman2011review}, maintenance planning~\cite{marwali1997integrated}, and transmission expansion~\cite{lumbreras2016new}. %
It is known that using DC power flow in transmission switching may produce infeasible or sub-optimal solutions, and some methods correct DC models for AC-infeasiblility~\cite{barrows2014correcting} or use convex relaxations such as Second-Order-Cone (SOC) 
\cite{bai2016two} to address this limitation. %
It is worth noting that existing methods for AC-feasibility recovery from DC-based solutions or heuristic solutions to transmission switching~\cite{crozier2022feasible, johnson2020k, capitanescu2014ac} are not directly applicable to OPS problems because OPS problems 1) tend to turn off many more power lines, 2) may involve load shed, and 3) may split the system into multiple islands.

Other problems are more closely related to the OPS problem 
in that large regions of the grid may be de-energized, the network may contain islands, and load shed may be required to find feasible power flow solutions. Examples 
include the \emph{Maximal Load Delivery} (MLD) problem in a severely damaged power system~\cite{mld}, where some lines may have to be disconnected for feasibility, as well as the associated problem of post-disaster restoration~\cite{rhodes2021powermodelsrestoration, rhodes2022role}, which decides how a set of damaged lines should be prioritized for repair. 
Both~\cite{rhodes2021powermodelsrestoration, rhodes2022role} and ~\cite{mld} include investigations into the choice of power flow formulations. They demonstrate that DC power flow may lead to sub-optimal binary decisions due to lack of consideration of reactive power and voltage problems, and that AC-based formulations may get stuck in local optima. SOC-based formulations seem to provide the best accuracy-quality trade-off given current solver technology, but can still be prohibitively slow. The results in~\cite{rhodes2021powermodelsrestoration} also demonstrated the need to redispatch generation with an AC-based optimal power flow after the binary decisions are made.

In this paper, we leverage the experiences from prior work to study the solution accuracy, ability to recover an AC feasible solution and solution speed in the context of the OPS problem. This problem differs from MLD and restoration problems in that we choose which power lines to de-energize from all lines in the system (as opposed to outage scenarios where a limited number of outaged lines are given) and that it contains a different objective with multiple parts (i.e. to balance wildfire risk reductions with load shed). Furthermore, prior work~\cite{rhodes2023security} has shown that the partially de-energized grids tend to have a primarily radial structure. These differences may have significant impact on the solution speed and accuracy of approximations at the optimal solution, and thus warrants an independent evaluation. 

In summary, the contributions of the work are as follows.  
1) We develop OPS formulations that use the AC, SOC, Network Flow (NF), and DC power flow formulations. \edit{Given the de-energization decisions obtained by solving OPS with different (non-AC) power flow formulations, we use an AC-Redispatch problem to recover AC-feasible power flow solutions.} %
2) We perform a comprehensive analysis of solution time and solution quality using 11 test networks with hundreds of risk scenarios. We assess the objective values and accuracy of solutions obtained with different power flow formulations, as well as the solution time.

The remainder of the paper is organized as follows: Section \ref{sec:formulation} introduces the OPS problem formulations and Section \ref{sec:ac_feasible} introduce the AC-feasible power flow recovery.  Section \ref{sec:case_study} analyzes the impact of power flow formulation on OPS, and Section \ref{sec:conclusion} concludes the work.

\section{Optimal Power Shutoff Problem Formulation} \label{sec:formulation}
In this section we introduce %
the OPS problem with AC, SOC, DC, and NF power flow formulations. The formulation is based on the DC-OPS problem that was first presented in~\cite{rhodes2020balancing}, while the adaptation to AC, SOC and NF formulations are new.

Across formulations, we consider a network with buses $ i \in \mathcal{B}$, generators $ g \in \mathcal{G}$, lines $ij \in \mathcal{L}$, load $d \in \mathcal{D}$, and shunts $s \in \mathcal{S}$.  Subsets of components (such as generators) at bus $i$ are defined as $\mathcal{B}^\mathcal{G}_i$. Parameters are \textbf{bold}  while variables are non-bold.  Binary variables representing the energization state of a component are $z\in\{0,1\}$ and indexed according to the component. If $z\!\!=\!\!1$, the component is energized, otherwise $z\!=\!0$. 

\subsection{AC Optimal Power Shutoff}
We start by defining the OPS problem with AC power flow equations, then show how the formulation must adapt for the other power flow formulations.

\subsubsection{Objective function}
The objective of the OPS problem 
seeks to maximize the power delivered while minimizing the risk of a wildfire ignition, with a pre-specified weighting factor $\boldsymbol{\alpha}\in[0,1]$ to determine a trade-off between load delivery and risk reduction. The objective function is given by
\begin{equation}
    \max \;\; (1-\boldsymbol{{\alpha}})\frac{\sum_{d\in\mathcal{D}} x_d \boldsymbol{w}_d \boldsymbol{P}^D_d}{\boldsymbol{P}^D_{tot}} - \boldsymbol{{\alpha}}\frac{\sum_{ij\in\mathcal{L}}z_{ij} \boldsymbol{R}_{ij}}{\boldsymbol{R}_{tot}}  \label{eq:obj}
\end{equation}
The first term represents the total power delivered in the system. The power delivered to each node $d$ is given by the power demand $\boldsymbol{P}_d^D$ multiplied by the continuous variable $x_d\in[0,1]$.  
The weighting factor $\boldsymbol{w}_d$ can be used to prioritize loads such as emergency services and community centers. 

The second term represents the risk of wildfire ignitions. The risk associated with an ignition from line $ij$ is given by  parameter $\boldsymbol{R}_{ij}$, which is multiplied by the binary energization status of a power line $z_{ij}$. 
\edit{The risk parameter $\boldsymbol{R}_{ij}$ represents the wildfire risk in the area around the power line, reflecting external factors such as vegetation, geography, and weather. This underlying wildfire risk represents the \emph{potential} for large and damaging wildfires to occur, and is not impacted by the power system operations. What is impacted by power system operations is the probability of ignitions caused by power lines. Specifically, if a line is de-energized, i.e. $z_{ij}=0$, the risk of a wildfire ignition is zero as a de-energized power line cannot ignite a wildfire. Conversely, if the line remains energized, i.e. $z_{ij}=1$, the line has a non-zero risk of igniting a fire.}
If the line remains energized, $z_{ij}=1$ and the line has a non-zero risk.  The total wildfire ignition risk is the summation over all power lines.  
The two terms are normalized by the total load demand $\boldsymbol{P}^D_{tot}=\sum_{d\in\mathcal{D}}\boldsymbol{P}^D_d$ and total wildfire risk $\boldsymbol{R}_{tot}=\sum_{ij\in\mathcal{L}}\boldsymbol{R}_{ij}$ in the system, respectively.

\subsubsection{Energization constraints}
\edit{Most power line de-energizations are aimed at reducing wildfire risk. However, additional components may have to be de-energized to enable a feasible power flow solution, and as a result we include de-energization decisions for all buses, loads and generators as well.}
In some cases, the energization status of a component are constrained by the energization status of the components they are connected to.  
For example, if a bus $i$ is de-energized $z_i=0$, any connected loads and shunts have to be shed (i.e. $x_d=0$ \edit{and $x_s=0$}), and all generators and lines have to be de-energized (i.e. $z_g=0$ and $z_{ij}=0$). These relationships are described by the following constraints, 
\begin{subequations}
\begin{align}
        &  z_{ig} \le z_{i}  &&  \forall g \in \mathcal{B}^\mathcal{G}_i, \; \forall i \in \mathcal{B}   \label{eq:ops_gen_active}\\
        &  z_{i j } \le z_{i} &&  \forall ij \in \mathcal{B}_i^{\mathcal{L}}, \; \forall i \in \mathcal{B}   \label{eq:ops_line_active}\\
        &  x_{d} \le z_{i}  &&   \forall d \in \mathcal{B}^\mathcal{D}_i, \; \forall i \in \mathcal{B}  \label{eq:ops_load_active} \\
        &  x_{s} \le z_{i}  &&   \forall s \in \mathcal{B}^\mathcal{S}_i, \; \forall i \in \mathcal{B}  \label{eq:ops_shunt_active}
\end{align}
\label{eq:ops_relationships} \vspace{-1.0em}
\end{subequations}

\edit{As stated previously each of these variables vary from fully de-energized $0$ to fully-energized $1$.  The $z$ variables are binary, while $x$ variables are continuous allowing continuous load and shunt shed. These variables bounds are shown in eq. \eqref{eq:z_bounds}}
\begin{subequations}
    \begin{align}
        & z_g \in \{0,1\} & \forall g \in \mathcal{G} \label{eq:zg_bound} \\
        & z_{ij} \in \{0,1\} & \forall ij \in \mathcal{L} \label{eq:zij_bound} \\
        & z_i \in \{0,1\} & \forall i \in \mathcal{B} \label{eq:zi_bound} \\
        & 0 \le x_d \le 1 & \forall d \in \mathcal{D} \label{eq:xd_bound} \\
        & 0 \le x_s \le 1 & \forall s \in \mathcal{S} \label{eq:xs_bound} 
    \end{align}
    \label{eq:z_bounds}
\end{subequations}

\subsubsection{Generation constraints}
The active power output of a generator $P_g^G$ is constrained between is upper $\overline{\boldsymbol{P}^G_g}$ and lower $\underline{\boldsymbol{P}^G_g}$ power limits when the generator is energized, and constrained to $0$ when the generator is de-energized, shown in eq. \eqref{eq:ops_active_gen_limits}.  The reactive power $Q_g^G$ is similarly constrained in eq. \eqref{eq:ops_reactive_gen_limits}.  
\begin{equation}
        z_{g}\underline{\boldsymbol{P}^G_g} \le P_{g}^G \le z_{g} \overline{\boldsymbol{P}^G_g} \quad \forall g \in \mathcal{G} \label{eq:ops_active_gen_limits}
\end{equation}
\begin{equation}
        z_{g}\underline{\boldsymbol{Q}^G_g} \le Q_{g}^G \le z_{g} \overline{\boldsymbol{Q}^G_g} \quad \forall g \in \mathcal{G} \label{eq:ops_reactive_gen_limits}
\end{equation}
\subsubsection{Power flow constraints} Active and reactive power balance at each node are given by
\begin{equation}
        \sum_{g\in\mathcal{B}_i^\mathcal{G}}P_{g}^G - \!\!\!\!\!\! \sum_{(i, j)\in\mathcal{B}_i^\mathcal{L}} \!\!\!\!\!  P_{i j}^L - \!\!\!  \sum_{d\in\mathcal{B}_i^\mathcal{D}} \!\!\! x_{d} \boldsymbol{P}^D_d - \!\!\!
        \sum_{s \in \mathcal{B}_i^{\mathcal{S}}} \!\! \boldsymbol{g}_s V_i^2 x_s = 0 \; \forall i \in \mathcal{B}  \label{eq:ops_active_power_balance} 
\end{equation}
\begin{equation}
        \sum_{g\in\mathcal{B}_i^\mathcal{G}}Q_{g}^G - \!\!\!\!\!\! \sum_{(i, j)\in\mathcal{B}_i^\mathcal{L}} \!\!\!\!\!  Q_{i j}^L - \!\!\!  \sum_{d\in\mathcal{B}_i^\mathcal{D}} \!\!\! x_{d} \boldsymbol{Q}^D_d + \!\!\!
        \sum_{s \in \mathcal{B}_i^{\mathcal{S}}} \!\! \boldsymbol{b}_s V_i^2 x_s = 0 \; \forall i \in \mathcal{B}  \label{eq:ops_reactive_power_balance} 
\end{equation}
where $P_{ij}^L, Q_{ij}^L$ represent the active and reactive power flow on line $ij$ and $\boldsymbol{Q}_d^D$ is the reactive power demand. \edit{ The shunt conductance and susceptance is given by $\boldsymbol{g}_s$ and $\boldsymbol{b}_s$ respectively, and are multiplied by the square of the node voltage magnitude $V_i$ and the shunt shed variable $x_s$}.

The thermal power limit for a power line $\boldsymbol{T}_{i j}$ is expressed in terms of the squared magnitude of complex power, which must be less than the squared thermal limit of the power line when energized $z_{ij}=1$, or set to $0$ when the line is de-energized $z_{ij}=0$. \edit{This limit is applied to power flow in both directions,}
\begin{subequations}
\begin{align}
        && 0 \le  (P_{i j}^L)^2 + (Q_{i j}^L)^2 \le \boldsymbol{T}_{i j}^2z_{i j} \quad \forall ij \in \mathcal{L} \label{eq:ops_complex_thermal_limit_to} \\
        && 0 \le  (P_{j i}^L)^2 + (Q_{j i}^L)^2 \le \boldsymbol{T}_{i j}^2z_{i j} \quad \forall ij \in \mathcal{L} \label{eq:ops_complex_thermal_limit_fr}
\end{align} \label{eq:ops_complex_thermal_limit}
\end{subequations}
The bus voltage magnitude $V_i$ must be within the upper $\boldsymbol{\overline{V}}_i$ and lower $\boldsymbol{\underline{V}}_i$ limits when bus $i$ is energized, and $0$ otherwise, 
\begin{equation}
    z_i \boldsymbol{\underline{V}}_i \le V_i \le z_i \boldsymbol{\overline{V}}_i \quad \forall i \in \mathcal{B} \label{eq:ops_voltage_magnitude}
\end{equation}

The AC power flow equations with line switching are
\begin{subequations}
\begin{align}
    \begin{split}
    P_{ij}^L &= z_{ij}  \Bigg( \Bigg. 
        \frac{\boldsymbol{g}_{ij} + \boldsymbol{g}_{i}}{|\boldsymbol{t}_{ij}|^2} V_i^2 %
        \!\! + \! \frac{-\boldsymbol{g}_{ij} \boldsymbol{t}^R_{ij} \!+\! \boldsymbol{b}_{ij} \boldsymbol{t}^I_{ij}}{|\boldsymbol{t}_{ij}|^2} V_{i} V_{j} \cos\left(\theta_{i} \!-\! \theta_{j}\right) \\ 
        & \!\!\!\! + \! \frac{-\boldsymbol{b}_{ij} \boldsymbol{t}^R_{ij} \!-\! \boldsymbol{g}_{ij} \boldsymbol{t}^I_{ij}}{|\boldsymbol{t}_{ij}|^2} V_{i}  V_{j} \sin\left(\theta_{i} \!-\! \theta_{j} \right) 
    \!\! \Bigg. \Bigg) \: \forall ij \in \mathcal{L} \label{eq:ac_p_power_flow_fr}
    \end{split}\\
    \begin{split}
    P_{ji}^L &= z_{ij}  \Bigg( \Bigg. 
        \left( \boldsymbol{g}_{ij} + \boldsymbol{g}_{j}\right) V_j^2 
        \!\! + \! \frac{-\boldsymbol{g}_{ij} \boldsymbol{t}^R_{ij} \!-\! \boldsymbol{b}_{ij} \boldsymbol{t}^I_{ij}}{|\boldsymbol{t}_{ij}|^2} V_{j} V_{i} \cos\left(\theta_{j} \!-\! \theta_{i}\right) \\
        & \!\!\!\! + \! \frac{-\boldsymbol{b}_{ij} \boldsymbol{t}^R_{ij} \!+\! \boldsymbol{g}_{ij} \boldsymbol{t}^I_{ij}}{|\boldsymbol{t}_{ij}|^2} V_{j}  V_{i} \sin\left(\theta_{j} \!-\! \theta_{i} \right) 
    \!\! \Bigg. \Bigg) \: \forall ij \in \mathcal{L}  \label{eq:ac_p_power_flow_to}
    \end{split}\\
    \begin{split}
    Q_{ij}^L &= z_{ij}  \Bigg( \Bigg.
         -\frac{\boldsymbol{b}_{ij} + \boldsymbol{b}_{i}}{|\boldsymbol{t}_{ij}|^2} V_i^2 
         \!\! - \! \frac{-\boldsymbol{b}_{ij} \boldsymbol{t}^R_{ij} \!-\! \boldsymbol{g}_{ij} \boldsymbol{t}^I_{ij}}{|\boldsymbol{t}_{ij}|^2} V_{i} V_{j} \cos\left(\theta_{i} \!-\! \theta_{j}\right) \\ 
        & \!\!\!\! + \! \frac{-\boldsymbol{g}_{ij} \boldsymbol{t}^R_{ij} \!+\! \boldsymbol{b}_{ij} \boldsymbol{t}^I_{ij}}{|\boldsymbol{t}_{ij}|^2} V_{i} V_{j} \sin\left(\theta_{i} \!-\! \theta_{j}\right) 
    \!\! \Bigg. \Bigg) \: \forall ij \in \mathcal{L} \label{eq:ac_q_power_flow_fr} 
    \end{split}\\
    \begin{split}
    Q_{ji}^L &= z_{ij}  \Bigg( \Bigg.
         - \left( \boldsymbol{b}_{ij} + \boldsymbol{b}_{j} \right) V_j^2 
         \!\! - \! \frac{-\boldsymbol{b}_{ij} \boldsymbol{t}^R_{ij} \!+\! \boldsymbol{g}_{ij} \boldsymbol{t}^I_{ij}}{|\boldsymbol{t}_{ij}|^2} V_{j} V_{i} \cos\left(\theta_{j}\!-\!\theta_{i}\right) \\ 
        & \!\!\!\! + \! \frac{-\boldsymbol{g}_{ij} \boldsymbol{t}^R_{ij} \!-\! \boldsymbol{b}_{ij} \boldsymbol{t}^I_{ij}}{|\boldsymbol{t}_{ij}|^2} V_{j} V_{i} \sin\left(\theta_{j}\!-\!\theta_{i}\right) 
    \!\! \Bigg. \Bigg) \: \forall ij \in \mathcal{L}  \label{eq:ac_q_power_flow_to}
    \end{split}
\end{align}
\label{eq:ops_ac_powerflow}
\end{subequations}
Parameters $\boldsymbol{g}_{ij}$ and $\boldsymbol{b}_{ij}$ are the conductance and suseptance, and $V_i, V_j$ and $\theta_i, \theta_j$ are the voltage magnitudes and angles at either end of the power line. \edit{Transformers are lossless and located at the $i$ side of the line with a fixed, complex-value voltage transformation $\boldsymbol{t}_{ij}$.  The real component of the transformation $\Re (\boldsymbol{t}_{ij})$ is $\boldsymbol{t}^R_{ij}$,  and the imaginary component of the transformation $\Im (\boldsymbol{t}_{ij})$ is $\boldsymbol{t}^I_{ij}$.}
When a power line is energized $z_{ij}=1$ the equations represent ordinary  AC power flow, while the power flow across the line is constrained to $0$ when the line is de-energized $z_{ij}=0$.

The voltage angle difference between two ends of a line must be between the maximum angle difference limits of the power line \edit{$\overline{\boldsymbol{\theta}_{ij}}$} when the line is energized, and should be unconstrained when the line is de-energized. This is expressed by the following big-M constraint,
\begin{subequations}
\begin{align}
    & \theta_{i} - \theta_{j} \le \overline{\boldsymbol{\theta}_{ij}}+ \boldsymbol{\theta}^{\Delta}_{max} (1-z_{i j}) \label{eq:ops_votlage_angle_max} \\
    & \theta_{i} - \theta_{j} \ge -\overline{\boldsymbol{\theta}_{ij}}- \boldsymbol{\theta}^{\Delta}_{max} (1-z_{i j}) \label{eq:ops_votlage_angle_min}
\end{align} \label{eq:ops_voltage_angle}
\end{subequations}
Here, $\boldsymbol{\theta}^{\Delta}_{max}$ is pre-computed to be the maximum angle difference between any two buses in the power system and used as the big-M value. %
The above equations are used to define the AC Optimal Power Shutoff (AC-OPS) problem, i.e. 
\begin{align}
    &\max &&  \mbox{Objective \eqref{eq:obj}} \nonumber\\
&\mbox{s.t.: \,\,\,}  \tag{AC-OPS}
&& \mbox{Component relationships:}~\eqref{eq:ops_relationships}, \eqref{eq:z_bounds} \nonumber \\[-2pt]
&&& \mbox{Generation constraints:}~\eqref{eq:ops_active_gen_limits}, \eqref{eq:ops_reactive_gen_limits}  \nonumber\\[-2pt]
&&& \mbox{AC power flow constraints:}~\eqref{eq:ops_active_power_balance} - \eqref{eq:ops_voltage_angle} \nonumber
\end{align}
This is a mixed-integer non-linear (and non-convex) problem (MINLP), which is a challenging class of problems to solve.

\subsection{SOC Optimal Power Shutoff}
We next describe the SOC-OPS problem, which uses the SOC relaxation of AC power flow, such that the objective value of the formulation is an upper bound on the optimal solution. 
The SOC power flow equations with component switching are adopted from~\cite{mld}, \edit{including several implementation aspects (i.e. specific formulations of constraints and extra variables that improve the MIP solver performance) that are derived in \cite{coffrin2015qc}}. Here we focus on describing the variables and constraints that differ from the AC power flow formulation, and refer the reader to \cite{jabr2006radial} for details on the derivation of the SOC relaxation. To obtain a convex relaxation of the power flow equations, the SOC formulation represents products between voltages using a set of lifted voltage squared variables \edit{such that $V_i V^*_i = W_{ii}$  for bus voltage. Voltage products across a power line are $V_i V_j^* =  W_{ij}$, where the variables in the model are the real and imaginary components given by $W_{ij}^R$, $W_{ij}^I$.  Additional variables that are added to improve the performance of a MIP solver are $W^{Fr}_{ij}$ and $W^{To}_{ij}$ which represent the squared bus voltage on either side of a line.  The variables $W^{Fr}_{ij}$ and $W^{To}_{ij}$ are equal to $W_{ii}$ and $W_{jj}$ when a line is energized, but differs in the case of de-energization. %
A variable for the voltage at a shunt $W^S_s$ is also required to allow a convex relaxation of shedding shunt power.  Note that the formulation does not contain any voltage angle variables.}

The upper and lower bounds of the squared nodal voltages \edit{$W_{ii}$ are} 
\begin{equation}
        z_{i} \boldsymbol{\underline{V}}^2_i \le  W_{ii}  \le  z_{i} \boldsymbol{\overline{V}}^2_i \quad  \forall i \in \mathcal{B} \label{eq:ops_soc_voltage_bus}
\end{equation}
\edit{Similarly, the bounds of the $W^{Fr}_{ij}$ and $W^{To}_{ij}$ are}
\begin{subequations}
    \begin{align}
        & z_{ij} \boldsymbol{\underline{V}}^2_i \le  W^{Fr}_{ij}  \le  z_{ij} \boldsymbol{\overline{V}}^2_i & \forall ij \in \mathcal{L} \\
        & z_{ij} \boldsymbol{\underline{V}}^2_j \le  W^{To}_{ij}  \le  z_{ij} \boldsymbol{\overline{V}}^2_j & \forall ij \in \mathcal{L} 
    \end{align} \label{eq:ops_soc_voltage_fr_to}
\end{subequations}
\edit{The relationship between $W_{ii}$ and $W^{Fr}_{ij}$ (and between $W_{jj}$ and $W^{To}_{ij}$) is given in the following equations.} 
\begin{subequations}
    \begin{align}
        & W_{ii} \ge W^{Fr}_{ij}  \ge W_{ii} - \boldsymbol{\overline{V}}_i^2 (1-z_{ij}) & \forall ij \in \mathcal{L} \\
        & W_{jj} \ge W^{To}_{ij}  \ge W_{jj} - \boldsymbol{\overline{V}}_j^2 (1-z_{ij}) & \forall ij \in \mathcal{L} 
    \end{align}  \label{eq:ops_soc_voltage_fr_to_ii_jj}
\end{subequations}
\edit{This constraint requires that the variables are equal when a power line is energized, but sets the value of $W^{Fr}_{ij}$ to $0$ (in combination with eq. \eqref{eq:ops_soc_voltage_fr_to}) without restricting the bus voltage when the line is de-energized.  The same argument applies to $W^{To}_{ij}$ and $W_{jj}$.}
The voltage product $V_i V^*_j$ for each power line are given in rectangular form, with the real and imaginary parts having the following bounds. 
\begin{subequations}
    \begin{align}
& z_{ij}\boldsymbol{\underline{W}}^R_{ij} \le  W^R_{ij}  \le  z_{ij}\boldsymbol{\overline{W}}^R_{ij} \quad \forall ij \in \mathcal{L}  \\
& z_{ij}\boldsymbol{\overline{W}}^I_{ij} \le  W^I_{ij}  \le  z_{ij}\boldsymbol{\underline{W}}^I_{ij} \quad   \forall ij \in \mathcal{L}
        \end{align} \label{eq:ops_soc_ij_ri_bounds}
\end{subequations}
\edit{The value of the upper and lower limits depends on the voltage angle limit of the line. The calculation of these parameters is shown in the Appendix. }
\edit{Limits on the voltage angle difference introduces constraints on the relationship between $W^R_{ij}$ and $W^I_{ij}$, given by}
\begin{equation}
        \tan \left(\boldsymbol{\underline{\theta}}_{ij}\right) W^R_{ij} \le W^I_{ij} \le \tan\left(\boldsymbol{\overline{\theta}}_{ij} \right) W^R_{ij} \quad \forall ij \in \mathcal{L}
    \label{eq:ops_soc_angle_limits}
\end{equation}

\edit{The line voltage variables are linked with the nodal voltage variables, and their relationship depends on the energization state of the line. The relationship is described by the equation $\left(W^R_{ij}\right)^2 +  \left(W^I_{ij}\right)^2  = W_{ii} W_{jj} z_{ij} \; \forall ij \in \mathcal{L}$, however, this equation is nonconvex and introduces a cubic equation on the right hand side of the constraint. To ensure convexity, the constraint is first relaxed by replacing the equality with an inequality  $\left(W^R_{ij}\right)^2 +  \left(W^I_{ij}\right)^2  \le W_{ii} W_{jj} z_{ij} \; \forall ij \in \mathcal{L}$, before the cubic term is removed by introducing the following relaxed inequalities}
\begin{subequations}
    \begin{align}
        & \left(W^R_{ij}\right)^2 +  \left(W^I_{ij}\right)^2  \le W_{ii} W_{jj}  & \forall ij \in \mathcal{L} \\
        & \left(W^R_{ij}\right)^2 +  \left(W^I_{ij}\right)^2  \le W_{ii} \overline{\boldsymbol{W_{jj}}} z_{ij} & \forall ij \in \mathcal{L} \\
        & \left(W^R_{ij}\right)^2 +  \left(W^I_{ij}\right)^2  \le \overline{\boldsymbol{W_{ii}}} W_{jj} z_{ij}& \forall ij \in \mathcal{L} 
    \end{align} \label{eq:ops_soc_ii_jj_ri}
\end{subequations}
\edit{The voltage across a shunt element varies according the amount the shunt is shed $x_s$.  The exact equation $W^S_s = W_{ii} x_s$ is nonconvex, and instead a McCormick relaxation is used}
\begin{subequations}
    \begin{align}
        &  W^S_s \ge 0  \quad \forall s \in \mathcal{B}^\mathcal{S}_i, & \forall i \in \mathcal{B} \\
        &  W^S_s \ge \boldsymbol{\overline{V}}^2_i (x_s -1 ) + W_{ii}  & \forall s \in \mathcal{B}^\mathcal{S}_i, \; \forall i \in \mathcal{B} \\
        &  W^S_s \le W_{ii}  \quad \forall s \in \mathcal{B}^\mathcal{S}_i, & \forall i \in \mathcal{B} \\
        &  W^S_s \le \boldsymbol{\overline{V}}^2_i x_s  \quad \forall s \in \mathcal{B}^\mathcal{S}_i, & \forall i \in \mathcal{B}
    \end{align} \label{eq:ops_soc_shunt_relaxation}
\end{subequations}
\edit{In the SOC nodal power balance constraints, $V_i^2$ is replaced by $W^S_{s}$ to describe the voltage across shunt elements,}
\begin{equation}
        \sum_{g\in\mathcal{B}_i^\mathcal{G}}P_{g}^G - \!\!\!\!\!\! \sum_{(i, j)\in\mathcal{B}_i^\mathcal{L}} \!\!\!\!\!  P_{i j}^L - \!\!\!\!  \sum_{d\in\mathcal{B}_i^\mathcal{D}} \!\!\! x_{d} \boldsymbol{P}^D_d  - \boldsymbol{g}_i W^S_s = 0 \quad \forall i \in \mathcal{B}  \label{eq:ops_soc_active_power_balance} 
\end{equation}
\begin{equation}
        \sum_{g\in\mathcal{B}_i^\mathcal{G}}Q_{g}^G - \!\!\!\!\!\! \sum_{(i, j)\in\mathcal{B}_i^\mathcal{L}} \!\!\!\!\!  Q_{i j}^L - \!\!\!\!  \sum_{d\in\mathcal{B}_i^\mathcal{D}} \!\!\! x_{d} \boldsymbol{Q}^D_d + \boldsymbol{b}_i W^S_s = 0 \quad \forall i \in \mathcal{B}  \label{eq:ops_soc_reactive_power_balance} 
\end{equation}
\edit{The SOC power flow equations with line switching use the $W$ variables, and are given by the following equations,}
\begin{subequations}
\begin{align}
    \begin{split}
        P_{ij}^L &=   
            \frac{\boldsymbol{g}_{ij} + \boldsymbol{g}_{i}}{|\boldsymbol{t}_{ij}|^2} W^{Fr}_{ij} 
            + \! \frac{-\boldsymbol{g}_{ij} \boldsymbol{t}^R_{ij} \!+\! \boldsymbol{b}_{ij} \boldsymbol{t}^I_{ij}}{|\boldsymbol{t}_{ij}|^2} W^R_{ij} \\ 
            & 
            + \! \frac{-\boldsymbol{b}_{ij} \boldsymbol{t}^R_{ij} \!-\! \boldsymbol{g}_{ij} \boldsymbol{t}^I_{ij}}{|\boldsymbol{t}_{ij}|^2} W^I_{ij}
        \quad \forall ij \in \mathcal{L} \label{eq:soc_p_power_flow_fr}
        \end{split}\\
    \begin{split}
        P_{ji}^L &=    
            \left( \boldsymbol{g}_{ij} + \boldsymbol{g}_{j} \right) W^{To}_{ij} 
            + \! \frac{-\boldsymbol{g}_{ij} \boldsymbol{t}^R_{ij} \!-\! \boldsymbol{b}_{ij} \boldsymbol{t}^I_{ij}}{|\boldsymbol{t}_{ij}|^2} W^R_{ij} \\
            & 
            + \! \frac{-\boldsymbol{b}_{ij} \boldsymbol{t}^R_{ij} \!+\! \boldsymbol{g}_{ij} \boldsymbol{t}^I_{ij}}{|\boldsymbol{t}_{ij}|^2} W^I_{ij}
        \quad \forall ij \in \mathcal{L}  \label{eq:soc_p_power_flow_to}
        \end{split}\\
    \begin{split}
        Q_{ij}^L &=    
             -\frac{\boldsymbol{b}_{ij} + \boldsymbol{b}_{i}}{|\boldsymbol{t}_{ij}|^2} W^{Fr}_{ij} 
            - \! \frac{-\boldsymbol{b}_{ij} \boldsymbol{t}^R_{ij} \!-\! \boldsymbol{g}_{ij} \boldsymbol{t}^I_{ij}}{|\boldsymbol{t}_{ij}|^2} W^R_{ij} \\ 
            & 
            + \! \frac{-\boldsymbol{g}_{ij} \boldsymbol{t}^R_{ij} \!+\! \boldsymbol{b}_{ij} \boldsymbol{t}^I_{ij}}{|\boldsymbol{t}_{ij}|^2} W^I_{ij}
         \quad \forall ij \in \mathcal{L} \label{eq:soc_q_power_flow_fr} 
    \end{split}\\
    \begin{split}
        Q_{ji}^L &=    
             - \left( \boldsymbol{b}_{ij} + \boldsymbol{b}_{j} \right) W^{To}_{ij} 
            - \! \frac{-\boldsymbol{b}_{ij} \boldsymbol{t}^R_{ij} \!+\! \boldsymbol{g}_{ij} \boldsymbol{t}^I_{ij}}{|\boldsymbol{t}_{ij}|^2} W^R_{ij} \\ 
            & 
            + \! \frac{-\boldsymbol{g}_{ij} \boldsymbol{t}^R_{ij} \!-\! \boldsymbol{b}_{ij} \boldsymbol{t}^I_{ij}}{|\boldsymbol{t}_{ij}|^2} W^I_{ij}
         \quad \forall ij \in \mathcal{L}  \label{eq:soc_q_power_flow_to}
    \end{split}
\end{align}
\label{eq:ops_soc_powerflow}
\end{subequations}
\edit{Note that when the power line is de-energized $z_{ij}=0$, the voltage variables $W^{Fr}_{ij}$, $W^{To}_{ij}$, $W^R_{ij}$, and  $W^I_{ij}$ also go to zero, correctly resulting in no power flow while de-energized.}

With this, the full SOC-OPS formulation is 
\begin{align}
    &\max &&  \mbox{Objective \eqref{eq:obj}} \nonumber\\
&\mbox{s.t.: \,\,\,}  \tag{SOC-OPS}
&& \mbox{Component relationships:}~\eqref{eq:ops_relationships}, \eqref{eq:z_bounds} \nonumber \\[-2pt]
&&& \mbox{Generation constraints:}~\eqref{eq:ops_active_gen_limits}, \eqref{eq:ops_reactive_gen_limits}  \nonumber\\[-2pt]
&&& \mbox{SOC power flow:}~\eqref{eq:ops_complex_thermal_limit}, \eqref{eq:ops_soc_voltage_bus}-\eqref{eq:ops_soc_ij_ri_bounds}, \eqref{eq:ops_soc_angle_limits}-\eqref{eq:ops_soc_powerflow}  \nonumber
\nonumber
\end{align}
This is a mixed-integer second-order cone problem (MISOCP). It is convex and thus easier to solve than the AC-OPS, but is still a very challenging problem class.

\subsection{DC Optimal Power Shutoff}
The DC-OPS problem is the original version of the OPS problem first proposed in \cite{rhodes2020balancing}. It uses the DC power flow linear approximation to model the power flow, which only considers real power and does not model reactive power or voltage constraints.   
The power flow $P_{ij}^L$ is constrained be lower than the thermal limit $\boldsymbol{T}_{i j}$ when the line is energized, and to $0$ when the power line is de-energized, as described by 
\begin{equation}
        -\boldsymbol{T}_{i j}z_{i j} \le  P_{i j}^L \le \boldsymbol{T}_{i j}z_{i j}  \quad  \forall ij \in \mathcal{L} \label{eq:ops_active_thermal_limit}
\end{equation}
The DC power flow equations with line switching are %
\begin{subequations}
\begin{align}
        &  P_{i j}^L \le -\boldsymbol{b}_{i j} (\theta_{i} - \theta_{j} + \boldsymbol{\theta}^{\Delta}_{max} (1-z_{i j})) &\forall ij \in \mathcal{L} \label{eq:ops_dc_flow_limit1} \\
        &  P_{i j}^L \ge -\boldsymbol{b}_{i j} (\theta_{i} - \theta_{j} - \boldsymbol{\theta}^{\Delta}_{max} (1-z_{i j})) &\forall ij \in \mathcal{L},\label{eq:ops_dc_flow_limit2}
\end{align}
\label{eq:ops_dc_powerflow}
\end{subequations}
which simplify to the standard DC power flow when the line is energized $z_{ij}=1$, and ensures that the voltage angle difference $\theta_i - \theta_j$ is unconstrained when the line is de-energized by introducing the big-M value from eq. \eqref{eq:ops_voltage_angle}. 

\edit{The DC power balance equation uses a fixed voltage magnitude of $1$ p.u. for shunt element power.  We multiply the shunt power by $x_s$ in this formulation to model the shunt de-energization.}
\begin{equation}
        \sum_{g\in\mathcal{B}_i^\mathcal{G}}P_{g}^G - \!\!\!\!\!\! \sum_{(i, j)\in\mathcal{B}_i^\mathcal{L}} \!\!\!\!\!  P_{i j}^L - \!\!\!\!  \sum_{d\in\mathcal{B}_i^\mathcal{D}} \!\!\! x_{d} \boldsymbol{P}^D_d - \boldsymbol{g}_i  x_s = 0 \quad \forall i \in \mathcal{B}  \label{eq:ops_dc_active_power_balance} 
\end{equation}

The full DC OPS problem is given by %
\begin{align}
    &\max &&  \mbox{Objective \eqref{eq:obj}} \nonumber\\
&\mbox{s.t.: \,\,\,}  \tag{DC-OPS}
&& \mbox{Component relationships:}~\eqref{eq:ops_relationships}, \eqref{eq:z_bounds} \nonumber \\[-2pt]
&&& \mbox{Generation constraints:}~\eqref{eq:ops_active_gen_limits} \nonumber\\[-2pt]
&&& \mbox{DC power flow constraints:}~\eqref{eq:ops_voltage_angle},\eqref{eq:ops_active_thermal_limit},\eqref{eq:ops_dc_powerflow},\eqref{eq:ops_dc_active_power_balance} \nonumber
\end{align}
The DC-OPS is a mixed-integer linear problem (MILP) which is easier \edit{to solve} than the above two problems, but still challenging to solve for large power system cases.

\subsection{ NF Optimal Power Shutoff}
The Network Flow OPS problem replaces the physics based power flow equations with a network flow model \cite{coffrin2016network}. This model considers active power conservation on each node and assumes that power can be sent along any edge as long as the transmission capacity constraint \eqref{eq:ops_active_thermal_limit} is respected. \edit{We include this formulation to represent the most simplistic way to model the power network when accounting for wildfire risk.}  The resulting optimization model is given by
\begin{align}
    &\max &&  \mbox{\eqref{eq:obj}} \nonumber\\
&\mbox{s.t.: \,\,\,}  \tag{NF-OPS}
&& \mbox{Component relationships:}~\eqref{eq:ops_relationships}, \eqref{eq:z_bounds} \nonumber \\[-2pt]
&&& \mbox{Generation constraints:}~\eqref{eq:ops_active_gen_limits} \nonumber\\[-2pt]
&&& \mbox{Network flow constraints:}~\eqref{eq:ops_active_thermal_limit},\eqref{eq:ops_dc_active_power_balance}\nonumber
\end{align}
Similar to the DC-OPS, the NF-OPS is also a mixed-integer linear program (MILP), though with fewer constraints. 

\section{Recovering AC Feasible Solutions} \label{sec:ac_feasible}
The above OPS formulations determine a set of binary de-energization decisions for each line, as well as a set of continuous variables that describe the generation, load and power flow in the system. Since the solutions obtained with SOC, DC and NF formulations are typically not AC-feasible, the load shed predicted by the optimization problem for the given de-energization solution may be wrong. To evaluate the true load shed amount, we seek to identify the minimum AC-feasible load shed (or, conversely, the maximum AC-feasible load delivery) for the set of de-energization decisions determined by those formulations. 
Our approach to recover an AC feasible solution leverages the \emph{AC-Redispatch} problem~\cite{rhodes2021powermodelsrestoration} to find an AC-feasible power flow. %

The AC-Redispatch problem is a continuous AC optimal power flow problem which maximizes the load served, i.e.,%
\begin{equation}
    \max \qquad \sum_{d\in\mathcal{D}} x_d \boldsymbol{w}_d \boldsymbol{P}_d  \label{eq:load_objective}
\end{equation}
We solve this problem with the de-energization decisions fixed to the values $\boldsymbol{\hat z}$ from an OPS solution, 
\begin{subequations}
    \begin{align}
    & z_{ij} = \boldsymbol{\hat{z}}_{ij} &\forall ij \in \mathcal{L} \label{eq:redis_branch_fixed} \\ 
    & z_{i} = \boldsymbol{\hat{z}}_{i} &\forall i \in \mathcal{B} \label{eq:redis_bus_fixed}  \\ 
    & z_{g} = \boldsymbol{\hat{z}}_{g} &\forall g \in \mathcal{G} \label{eq:redis_gen_fixed}
    \end{align}
    \label{eq:redis_fixed}
\end{subequations}
For conciseness of notation, we use the same equations used in the AC-OPS problem, with added constraints to fix the binary values to the energization state from the solution of an OPS problem, shown in \eqref{eq:redis_fixed}. In practice, this problem is modeled as a continuous NLP problem rather than a MINLP problem.
The full AC-Redispatch problem is given by
\begin{align}
    &\max &&  \eqref{eq:load_objective} \nonumber \tag{AC-Redispatch}\\
&\mbox{s.t.: \,\,\,}
&& \mbox{Component relationships:}~\eqref{eq:redis_fixed}, \eqref{eq:xd_bound}, \eqref{eq:xs_bound} \nonumber \\[-2pt]
&&& \mbox{Generation constraints:}~\eqref{eq:ops_active_gen_limits}, \eqref{eq:ops_reactive_gen_limits}  \nonumber\\[-2pt]
&&& \mbox{AC Power flow constraints:}~\eqref{eq:ops_active_power_balance} - \eqref{eq:ops_voltage_angle} \nonumber
&&
\end{align}

\section{Case Study} \label{sec:case_study}
\begin{figure*}[t]
\vspace{-1.0em}
    \centering
    \begin{subfigure}{0.3\textwidth}
        \includegraphics[width=\columnwidth]{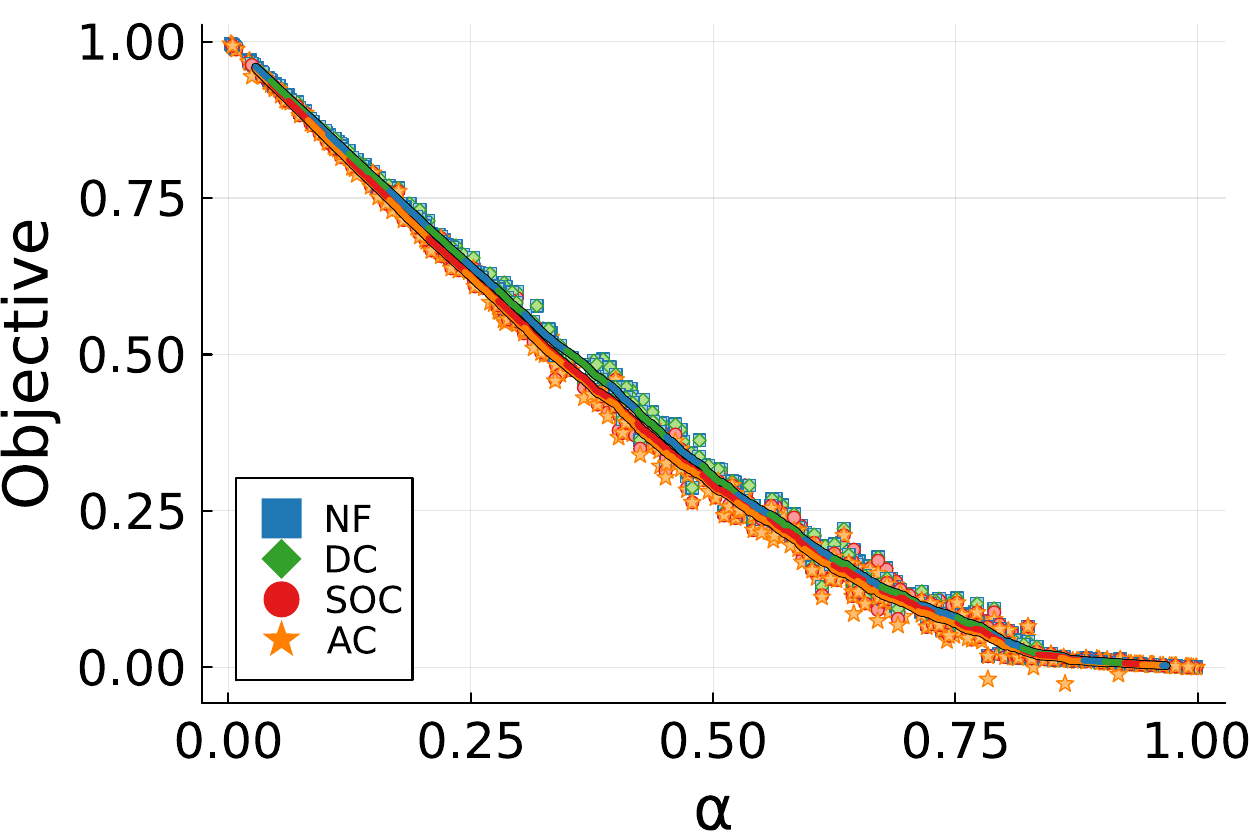}
        \caption{} \label{fig:case14_obj}
    \end{subfigure}   
    \begin{subfigure}{0.3\textwidth}
        \includegraphics[width=\columnwidth]{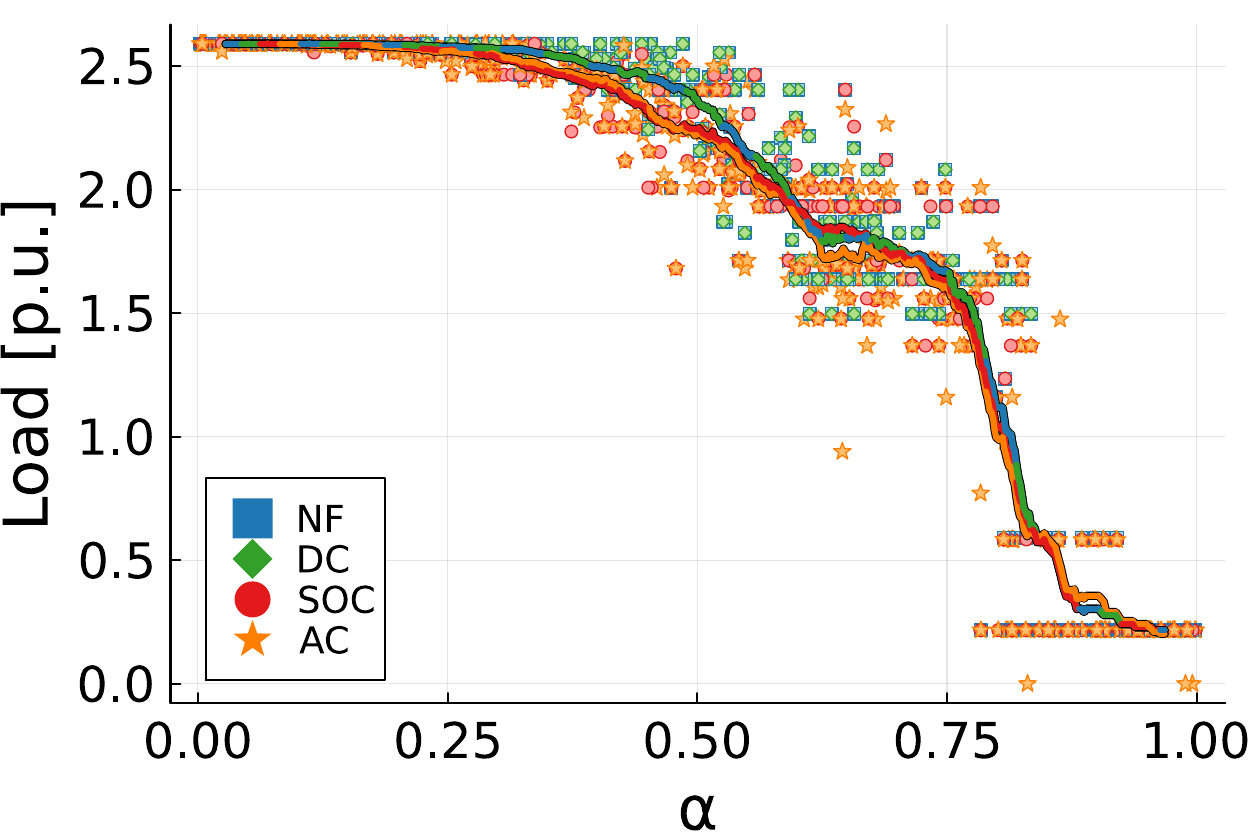}
        \caption{} \label{fig:case14_load}
    \end{subfigure} 
        \begin{subfigure}{0.3\textwidth}
        \includegraphics[width=\columnwidth]{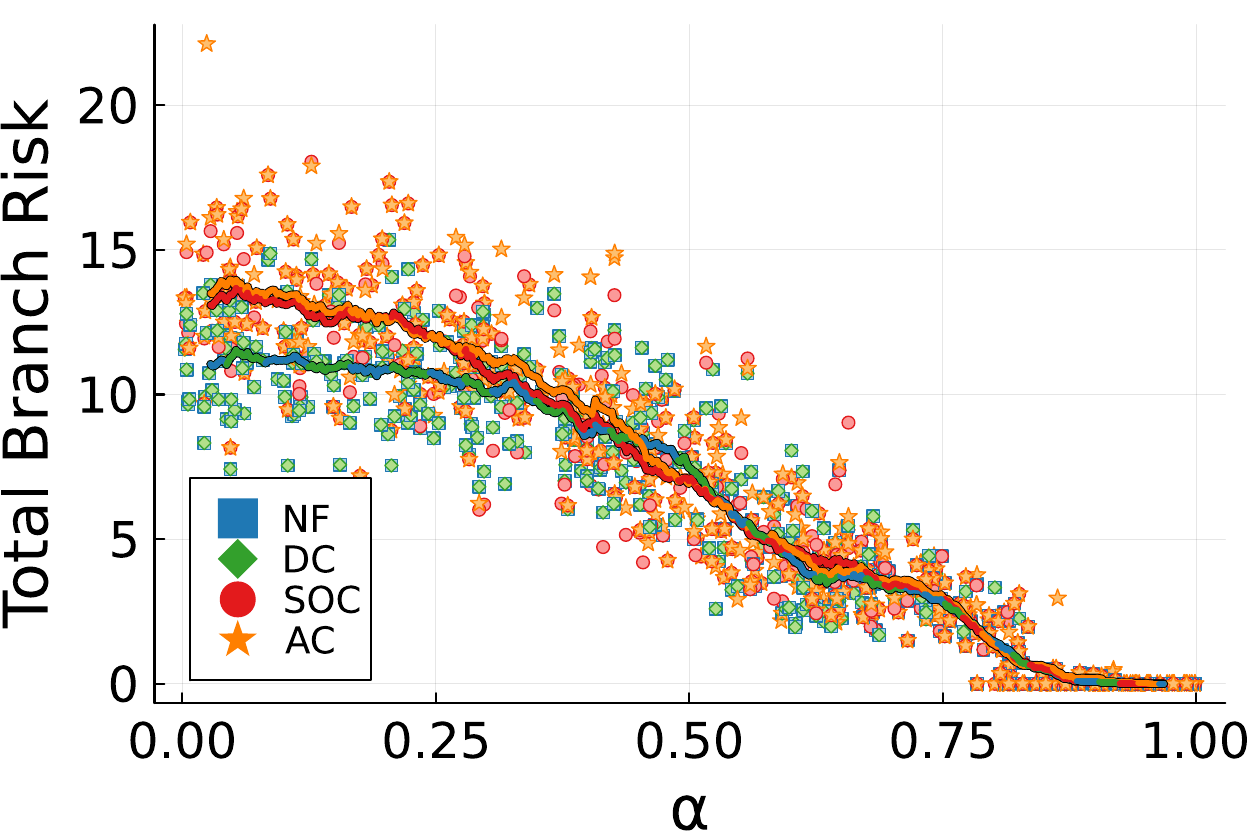}
        \caption{} \label{fig:case14_risk}
    \end{subfigure} 
    \begin{subfigure}{0.3\textwidth}
        \includegraphics[width=\columnwidth]{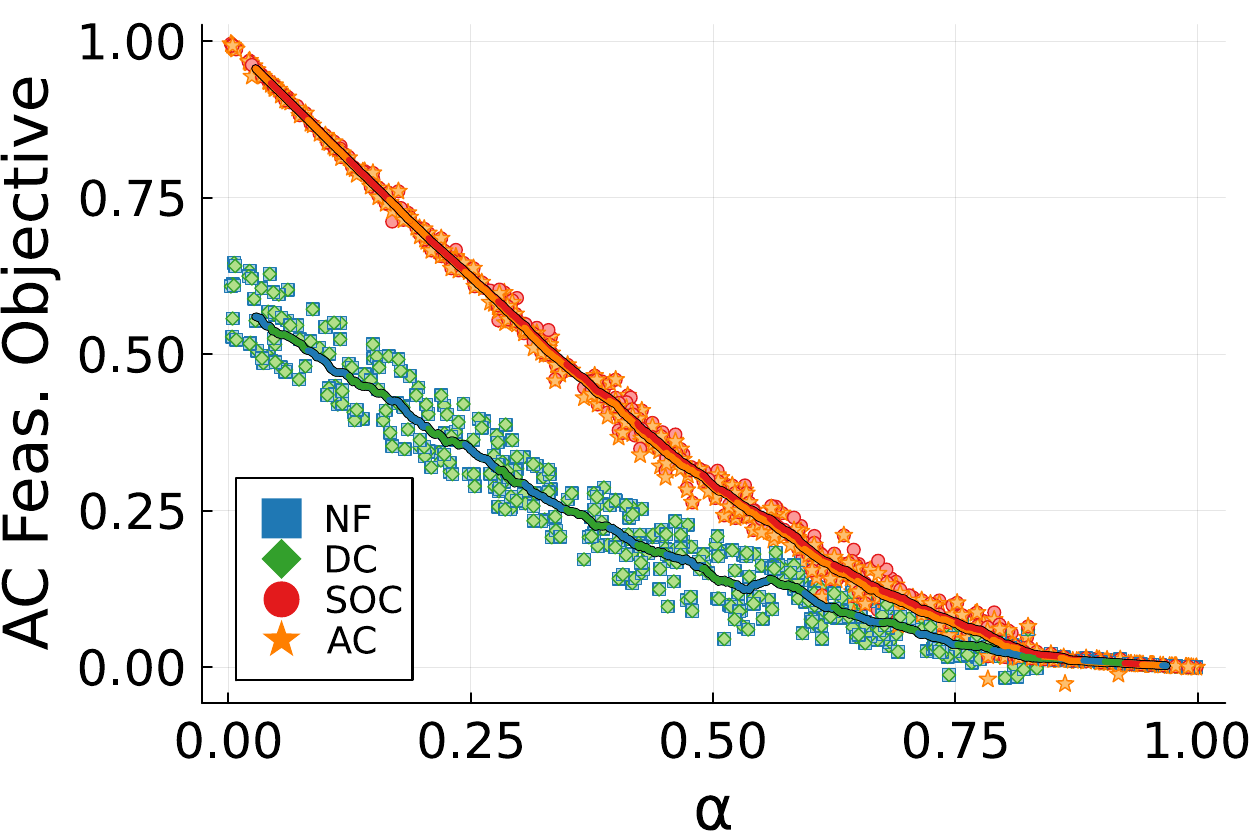}
        \caption{} \label{fig:case14_red_obj}
    \end{subfigure} 
    \begin{subfigure}{0.3\textwidth}
        \includegraphics[width=\columnwidth]{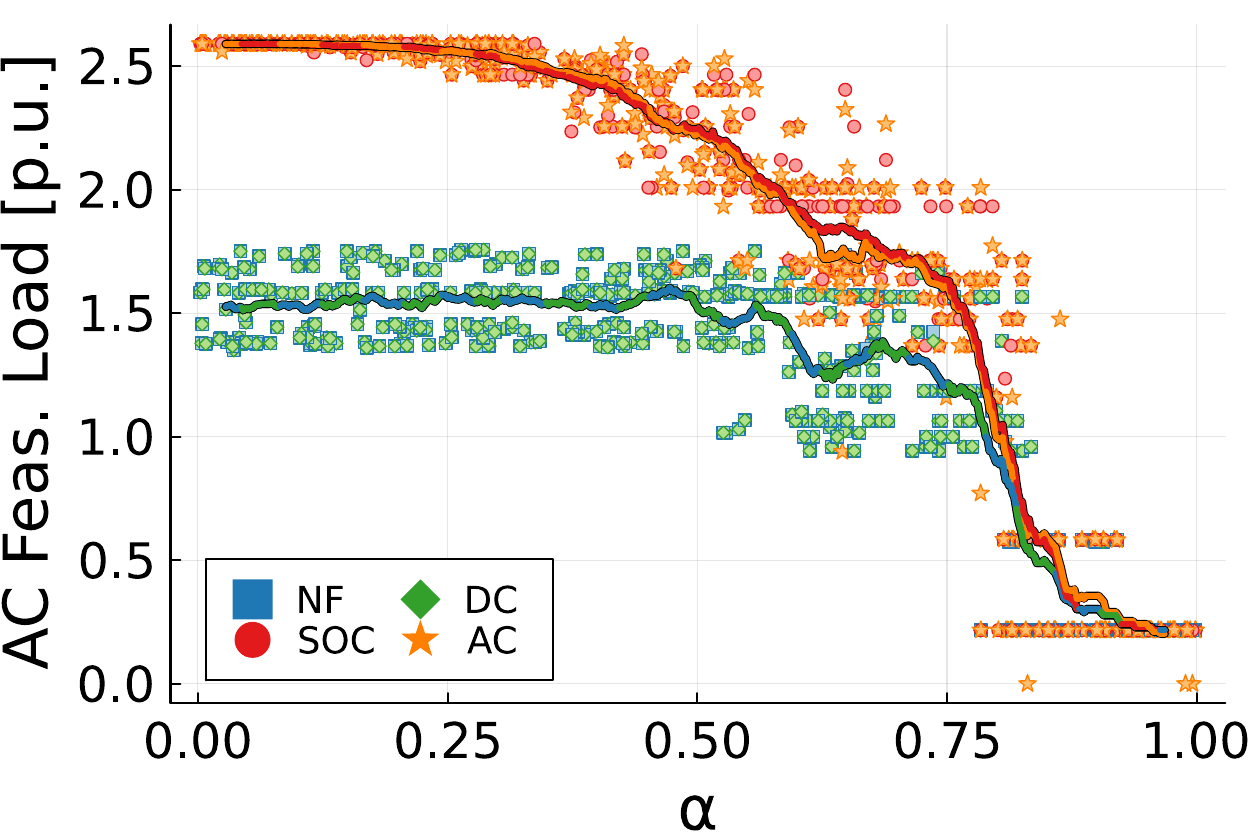}
        \caption{} \label{fig:case14_red_load}
    \end{subfigure}
    \begin{subfigure}{0.3\textwidth}
        \includegraphics[width=\columnwidth]{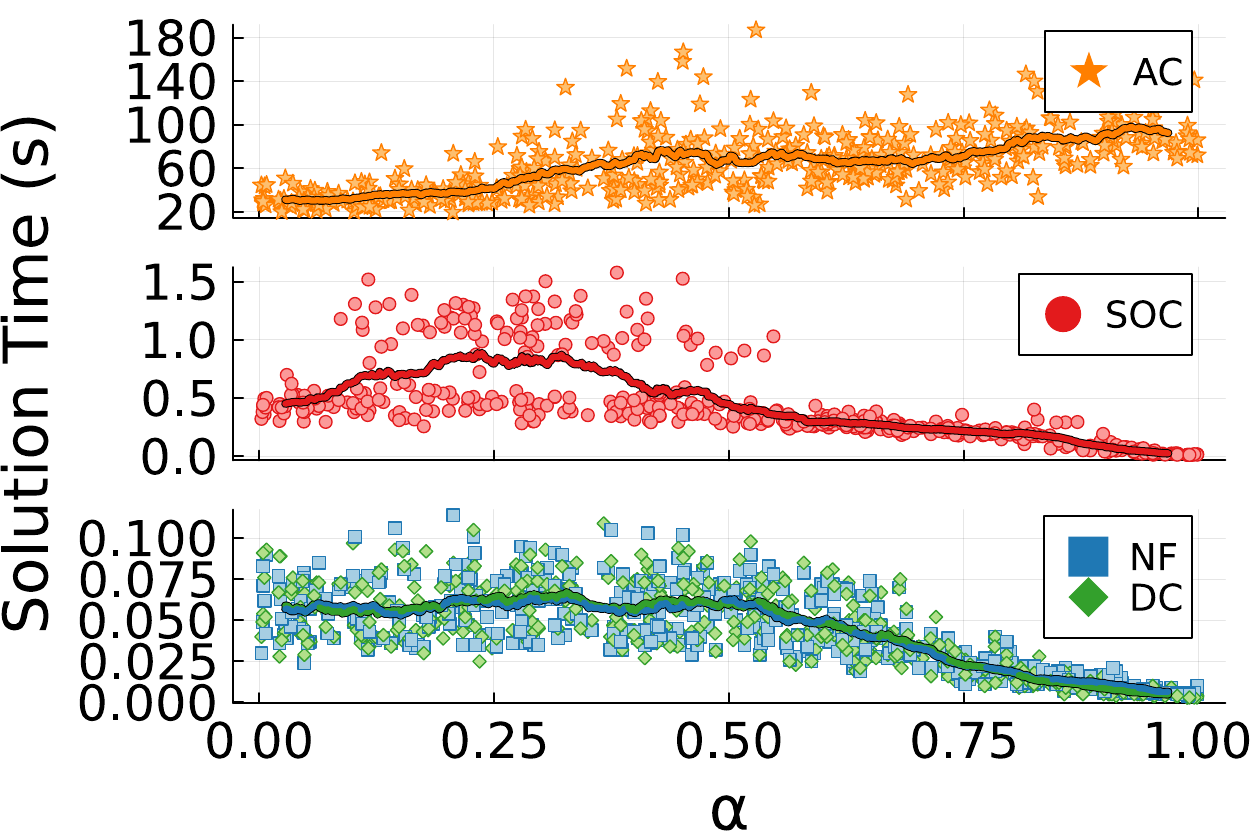}
        \caption{}
        \label{fig:case14_time}
    \end{subfigure}
    \caption{\small \textbf{Case IEEE-14} Fig. \ref{fig:case14_obj} shows a  scatterplot of the OPS objective (Load and Risk) for 500 risk scenarios, solved with four different power flow formulations, as a function of the $\alpha$ parameter. The power delivered to loads is shown in Fig. \ref{fig:case14_load} while the wildfire risk is shown in Fig. \ref{fig:case14_risk}.  Fig. \ref{fig:case14_red_obj} shows the re-calculated objective after an AC-feasible power flow is found.  The reduction in load served is shown in Fig. \ref{fig:case14_red_load}.  
    Fig. \ref{fig:case14_time} shows the solution time of the OPS problem for each formulation.
    }
    \label{fig:case14}
    \vspace{-1.5em}
\end{figure*}

We investigate the solution quality and solution time for the normalized OPS problem with different power flow formulations across a range of different power system test cases, wildfire values and choices of the trade-off parameter $\alpha$. 

\subsection{Case Study Set-Up}

\noindent \emph{1) Software Implementation} The OPS problem is implemented using the Julia programming language \cite{julia} with the JuMP optimization package \cite{jump}. We leverage the DC-OPS implementation in PowerModelsWildfire.jl \cite{rhodes2020balancing}, and extended the package to allow for AC, NF and SOC power flow formulations. These implementations are now publicly available. We use the solvers Gurobi \cite{gurobi} for the NF-, DC- and SOC-OPS problems, which are mixed-integer convex (linear or quadratic) programs, Juniper \cite{juniper} as a solver for the AC-OPS which is mixed-integer non-convex (non-linear) problems, and Ipopt \cite{ipopt} for the AC-Redispatch problem, which is a continuous non-linear problem. We use Distributions.jl \cite{Distributions} to generate random values for our input data (as further described below) and PowerPlots.jl \cite{powerplots} and Plots.jl \cite{plots} for visualization.

\noindent \emph{2) Test systems} We use 11 cases from PGLib \cite{pglib}, ranging from 3 to 118 buses. The names of the systems are listed along with the results in Table I. 

\noindent\emph{3) Wildfire data} 
The PGLib test cases do not include any detailed geographical information, and therefore %
there is no wildfire power line risk data available. To enable testing on a variety of systems and a large number of scenarios per system, we use randomly generated wildfire risk coefficients $\boldsymbol{R}_{ij}$ drawn from a Rayleigh distribution.  %
We chose a Rayleigh distribution because it is a good approximation of wind speed variation \cite{pishgar2015wind}, and local wildfire risk is closely correlated with wind speed \cite{rios2018interpolation}.

\subsection{Numerical Experiment Set-Up}
We generate 500 scenarios for each PGLib test case by sampling wildfire risk coefficients $\boldsymbol{R}_{ij}$ for each power line from the Rayleigh distribution. We also sample an $\alpha$ trade-off parameter from a uniform distribution for each scenario, which 
allows us to study the variation of these scenarios as alpha changes. \edit{In practice, the selection of $\alpha$ to weigh the trade off of load shed and wildfire ignition risk is a policy decision based on the priorities of the grid operator or regulator \cite{rhodes2020balancing}. We study a wide range of values of $\alpha$ to understand how it may impact solution quality in addition to its primary purpose of a trade off of priorities in the objective function.}  For each of these scenarios, we solve the NF-, DC-, SOC- and AC-OPS formulations.
To maintain a reasonable computational time when solving 500 scenarios, we enforce a 30 minute time limit for each optimization problem and include the results from the time limited scenarios. 
We do not include results for the AC-OPS on test cases larger than 14-buses and for the SOC-OPS on test cases larger than 39 buses, as a significant number of scenarios (more than 75 of 500) could not be solved within the time limit. 
For each solution, we record the
objective value,
load delivered,
total power line wildfire risk, 
and solution time. 

After solving the OPS problem, we use AC-Redispatch to recover an AC-feasible solution, and record the resulting load delivered. We note that the AC-Redispatch problem \edit{always found a locally optimal solution, however, that solution may be a trivial AC-feasible solution with no power flow on any lines.}

We first present detailed analysis on the IEEE 14 bus case before showing summarized results for the remaining systems. 

\subsection{IEEE 14 bus case} \label{sec:case14}
A collection of results for the IEEE 14-bus case are detailed in Fig. \ref{fig:case14}, where different values on the y-axis are plotted against the range of $\alpha$ values on the x-axis.
Individual solutions (corresponding to a set of risk coefficients and an $\alpha$ value) are shown as points in the scatter plots. The lines represent rolling averages of the nearest 30 data points to illustrate trends.

\subsubsection{Optimization Problem Results}
We first discuss the results obtained by solving the OPS.
The upper row shows the objective value in Fig. \ref{fig:case14_obj}, load delivered in Fig. \ref{fig:case14_load} and wildfire risk in Fig. \ref{fig:case14_risk}. 
We first observe from Fig. \ref{fig:case14_obj} that all power flow formulations find solutions with objective values in a similar range. However, the AC- and SOC-OPS solutions appear to have a slightly lower total objective value than the DC- and NF-OPS solutions. A closer look at Fig. \ref{fig:case14_load} and \ref{fig:case14_risk}, explain why. At $\alpha$ values $\!>\!\!0.5$, the AC- and SOC-OPS solutions serve slightly less load while the power line risk is slightly higher. At $\alpha$ values $\!<\!0.5$, the AC- and SOC-OPS solutions \edit{have similar power line risk and serve a similar amount of load as the DC- and NF-OPS solutions. This indicates that when conducting a small or moderate PSPS, the solutions obtained with AC- and SOC-power flow require a larger number of lines to be remain energized, thus finding solutions with a higher wildfire risk. At large-scale PSPS when most of the grid id de-energized, there is little difference between the linear formulations and the SOC or AC formulations, as the ability to delivery load is almost entirely due to the connectivity of the power grid.}

\subsubsection{AC-Feasible Load Delivery}
Next, we investigate how the load delivery changes after we recover an AC-feasible load delivery solution.
The results are shown in Fig. \ref{fig:case14_red_obj} and \ref{fig:case14_red_load}, which shows the re-calculated objective value and AC-feasible load delivery.
By comparing the AC-feasible load delivery in Fig. \ref{fig:case14_red_load} with the load delivery predicted by the optimization solution in Fig. \ref{fig:case14_load}, we observe that the load delivery achieved with the NF and DC solutions is significantly reduced for $\alpha\!<\!0.85$. The load delivery remains largely the same for the SOC solution across all $\alpha$ values (and the AC solution is already AC feasible). As a result, the objective values of the NF and DC-OPS solutions are significantly reduced, as seen in Fig. \ref{fig:case14_red_obj}.
This indicates that while the NF-, DC-OPS solutions initially produced similar (or even slightly better) objective values, the de-energization decisions were significantly different due to inaccuracies in the NF and DC power flow formulations. As a result, the NF- and DC-OPS solutions only deliver half as much load as the AC and SOC solutions once we recover an AC-feasible solution. This demonstrates that using a more detailed power flow model significantly improves the results. %

\subsubsection{Solution time}
Finally, Fig. \ref{fig:case14_time} shows the solution time needed to obtain the solution for each instance. Note that we plot the different problems on three separate scales. We observe that the NF- and DC-OPS problems have similar solution times, while the SOC-OPS take an order of magnitude, i.e., 10x longer to solve. The AC-OPS is slowest, taking 1000x  longer to solve than the DC and NF solutions.  Interestingly, the NF, DC, and SOC formulations are most challenging to solve when $\alpha \!<\! 0.5$  where solutions have little load shed and moderate risk reduction.  The solution time reduces significantly as $\alpha \rightarrow 1$ where the solution is total de-energization.  The trend is not the same for the AC-OPS problem, where the solution time increases as $\alpha \rightarrow 1$.

\subsection{PGLib Test Cases}\label{sec:all_cases}
\begin{table*}[t]
\vspace{-1.0em}
    \centering
    \caption{Average OPS Objective, AC-Feasible Objective, and Difference \edit{on 500 Scenarios}}
    \resizebox{\textwidth}{!}{%
    \begin{tabular}{ r|ccc|ccc|ccc|c} 
      {} & \multicolumn{3}{c|}{\textbf{NF}} & \multicolumn{3}{c|}{\textbf{DC}} & \multicolumn{3}{c|}{\textbf{SOC}} & \textbf{AC} \\
      \hline
      \textbf{Case} & \textbf{Obj.} & \textbf{AC-Feas. Obj.} & \textbf{Diff.} & \textbf{Obj.} & \textbf{AC-Feas. Obj.} & \textbf{Diff.} & \textbf{Obj.} & \textbf{AC-Feas. Obj.} & \textbf{Diff.} & \textbf{Obj.} \\ 
      \hline
      LMBD 3 Bus & .433293 & .413538 & .019755 & .433293 & .413538 & .019755 & .442606 & .442415 & .000191 & .436549 \\
      PJM 5 Bus & .391419 & .389826 & .001593 & .391419 & .389810 & .001610 & .391288 & .391088 & .000200 & .345657 \\
      IEEE 14 Bus & .369441 & .200461 & .168981 & .369441 & .200444 & .168998 & .359717 & .359597 & .000120 & .355732 \\
      IEEE RTS 24 Bus & .433100 & .367518 & .065582 & .433100 & .363706 & .069394 & .425365 & .424801 & .000564  & --\\
      AS 30 Bus & .351774 & .300583 & .051191 & .351774 & .300586 & .051188 & .348536 & .347950 & .000585 & --\\
      IEEE 30 Bus & .350683 & .194934 & .155749 & .350683 & .194945 & .155738 & .343981 & .342727 & .001253 & --\\
      EPRI 39 Bus & .348876 & .325593 & .023284 & .348876 & .325941 & .022935 & .345620 & .344390 & .001231 & --\\
      IEEE 57 Bus & .410583 & .374980 & .035603 & .410583 & .374980 & .035603 & -- & -- & -- & --\\
      IEEE RTS 73 Bus & .410938 & .311074 & .099864 & .410938 & .310687 & .100251 & -- & -- & -- & --\\
      PEGASE 89 Bus & .636862 & .503009 & .133853 & .636860 & .503633 & .133227 & -- & -- & -- & --\\
      IEEE 118 Bus & .342006 & .286755 & .055252 & .342007 & .286680 & .055326 & -- & -- & -- & --\\
    \end{tabular}
    }
    \label{tab:redis_diff}
\end{table*}

\begin{table*}
\parbox{.49\linewidth}{
\centering
    \caption{Difference of AC-Feasible Objectives}
        \resizebox{\linewidth}{!}{%
    \begin{tabular}{ r|cccc} 
      \textbf{Case} & \textbf{DC-NF} & \textbf{SOC-DC} & \textbf{AC-DC} & \textbf{SOC-AC} \\ 
      \hline
      LMBD 3 Bus & 0.0 & .028876 & .023011 & .005866 \\
      PJM 5 Bus & -.000016 & .001279 & -.044152 & .045431 \\
      IEEE 14 Bus & -.000017 & .159154 & .155288 & .003866 \\
      IEEE RTS 24 Bus & -.003812 & .061095 & -- & --\\
      AS 30 Bus & .000003 & .047364 & -- & --\\
      IEEE 30 Bus & .000011 & .147783 & -- & --\\
      EPRI 39 Bus & .000348 & .018449 & -- & --\\
      IEEE 57 Bus & -0.0 & -- & -- & --\\ 
      IEEE RTS 73 Bus & -.000387 & -- & -- & --\\
      PEGASE 89 Bus & .000623 & -- & -- & -- \\ 
      IEEE 118 Bus & -.000074 & -- & -- & --\\ 
    \end{tabular}
    }
    \label{tab:redis_compare}
}
\hfill
\parbox{.49\linewidth}{
\centering
    \caption{\edit{Number of Scenarios Where OPS Overestimates AC-Feasible Load Delivery By Greater Than $20\%$}}
    \resizebox{0.63\linewidth}{!}{%
    \begin{tabular}{ r|ccc}
      \textbf{Case} & \textbf{NF} & \textbf{DC} & \textbf{SOC} \\ 
      \hline
      LMBD 3 Bus & 0 & 0 & 0 \\
      PJM 5 Bus & 0 & 0 & 0 \\
      IEEE 14 Bus & 379 & 379 & 0 \\
      IEEE RTS 24 Bus & 77 & 79 & 2 \\ 
      AS 30 Bus & 36 & 36 & 0 \\
      IEEE 30 Bus & 366 & 366 & 0 \\
      EPRI 39 Bus & 17 & 17 & 1 \\ 
      IEEE 57 Bus & 1 & 1 & -- \\
      IEEE RTS 73 Bus & 160 & 168 & -- \\ 
      PEGASE 89 Bus & 144 & 144 & -- \\
      IEEE 118 Bus & 1 & 1 & -- \\ 
    \end{tabular}
    }
    \label{tab:infeasible} 
}
\vspace{-1.5em}
\end{table*}
We next show summary results for each of the PGLib test cases. 
Table \ref{tab:redis_diff} reports the summary metrics of solving the 500 scenarios on each of the 11 cases from PGLib using all formulations of the OPS problem.
For each network and formulation, we show the average objective value, the average objective value after finding an AC-feasible power flow, and the difference between the two. %
The number of scenarios where the AC-feasible solution has significant additional load-shed is shown in Table \ref{tab:infeasible}. %

\subsubsection{Accuracy of predicted load shed}
We first discuss the ability of the power flow formulations to accurately predict the load shed value, i.e. have a small difference between the optimization objective of the OPS problem and the objective value after the AC-feasibility recovery is solved.
For NF-OPS and DC-OPS, the difference between the objective value of the OPS problem and the re-evaluated objective using AC-feasible power flow range from \edit{$0.001$} on the PJM 5 bus network, to \edit{$0.168$} on IEEE 14 bus network. This indicates that on some networks, DC-OPF and NF-OPF are reasonably accurate on average, but on other networks they can serve less than half the estimated power delivery. SOC-OPS finds solutions that are very close to AC-feasible, with an average difference of less than \edit{0.001} for most networks.

\subsubsection{Quality of solutions} 
Next, we analyze the solution quality by comparing the amount of load shed after recovering an AC-feasible solution. 
Table \ref{tab:redis_compare} shows the average objective value difference of the formulations, after finding an AC-feasible power flow. We compare the solutions pairwise, 
with positive values indicating that the formulation listed first performs better. 
From the results, we first observe that the DC and NF objective values are nearly identical across all cases, with an average difference close to 0. We therefore compare only the DC formulation with the other power flow formulations. The difference between the DC and SOC objective values shows that the DC formulations performs worse that the SOC solution on all networks. Interestingly, the SOC-OPS with AC-feasibility recovery outperforms the AC-OPS problem on all cases, while even the DC-OPS performs better than the AC-OPS on the PJM 5 bus case. This is because the AC-OPS finds locally optimal shutoff solutions when using the Juniper solver.

\subsubsection{Recovering an AC-Feasible solution}

We next discuss the challenges involved with recovering AC-Feasible solutions. 
\edit{Table \ref{tab:infeasible} reports how many scenarios for each network had over a $20\%$ difference between the predicted OPS load shed and the AC-Redispatch load shed, meaning that if the planned PSPS was implemented it could result in over 20\% extra load shed for customers.}
We observe that the NF and DC formulations \edit{often significantly underestimate levels of load shed on most networks. For example, the IEEE 89 bus network contains more than 20\% extra load shed when solved with AC-Redispatch in 144 out of the 500 scenarios, while the IEEE 30 bus network has high additional load shed in 366 out of 500 scenarios.}

\edit{However, this is not the case on all networks. The solutions from the 57-bus and 118-bus networks only differ by greater than $20\%$ from the AC-feasible solution in 1 scenario out of 500 for each network. The  5-bus, and 3-bus networks never differ from the AC-Feasible solution by more than 20\%.  It is surprising that some of these networks can be modeled somewhat accurately with the linear formulations while other networks, while others have significant load estimate errors in over $70\%$ of the scenarios.}
\edit{It is however clear that the SOC formulation is more accurate than the DC or NF formulations. Across all seven networks,} the SOC formulation has only \edit{three scenarios (out of 3,500) where the AC-feasible solution was contained more than $20\%$ additional load shed.}

\if\removepowerplot1
\else
\begin{figure}[t]
\centering
\includegraphics[width=0.8\columnwidth]{fig/MLD Infeasible DC (case30_as).pdf}
\caption{\small \textbf{Case AS 30}: Lines de-energized by the DC-OPS problem are shown in grey.  Additional lines de-energized for AC-feasibility are orange.  Energized lines are red. }
\label{fig:MLD}
\vspace{-1.5em}
\end{figure}
To further illustrate how the DC and NF-OPS solutions fail, we show an example of using our proposed algorithm %
to recover an AC-feasible solution on the AS-30 network in Fig. \ref{fig:MLD}. \textcolor{blue}{XX} lines and \textcolor{blue}{YY}  generators in grey were de-energized by DC-OPS to reduce the wildfire risk, resulting in a radial network with very long chains of lines.  However, the resulting network did not have any AC-feasible power flow. The SOC-MLD problem was solved and de-energized an additional components  \textcolor{blue}{XX} lines and \textcolor{blue}{YY}  generators to find a feasible power flow. These de-energized components are shown in \textcolor{blue}{orange?}, and lead to several smaller radial systems with the final energized network configuration is shown in red. 
\fi

\subsubsection{Solution times}
\begin{figure}[t]
\vspace{-1.0em}
\includegraphics[width=\columnwidth]{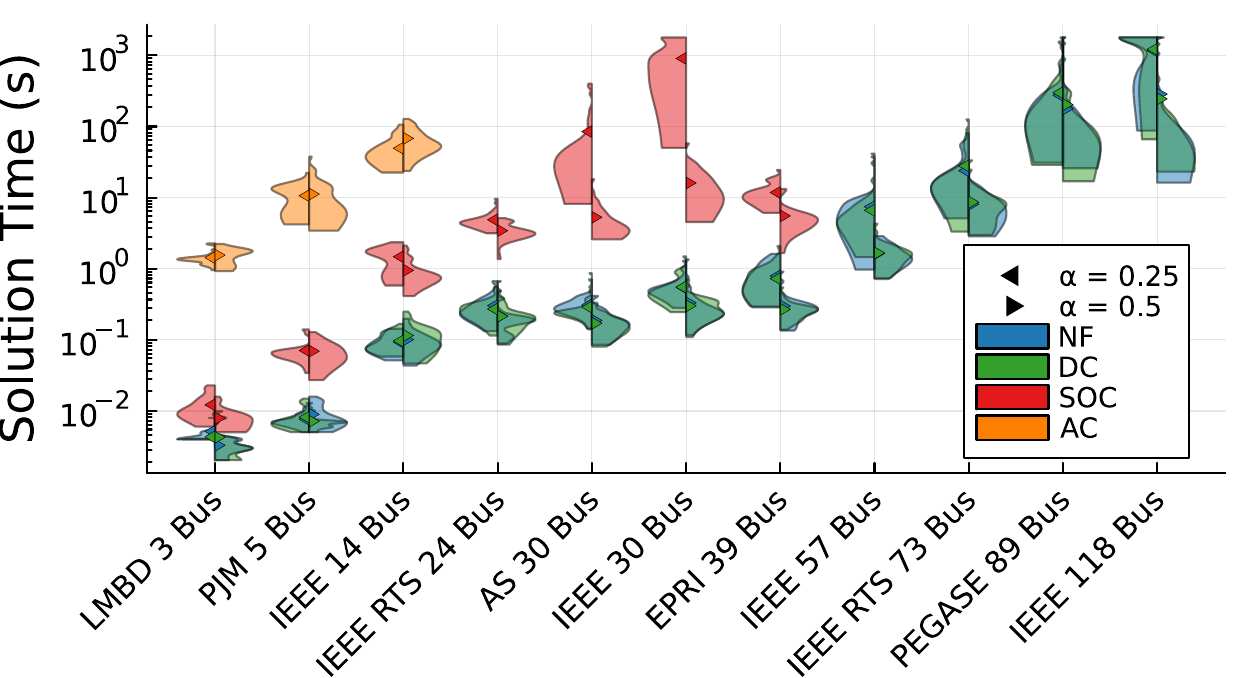}
\caption{\small \textbf{OPS Solve Speed}: Distribution of solution time for the PGLib cases at $\alpha$=0.25 and $\alpha$=0.5. The triangles denote the mean value.}
\label{fig:solve_times}
\vspace{-1.5em}
\end{figure}

Finally, we discuss the solution time. As seen for the IEEE 14 bus system, the solve time can vary greatly depending on the trade-off parameter $\alpha$ and the risk coefficients $\boldsymbol{R}_{ij}$.  
To evaluate solution speed for the PGLib test cases, we solve 50 scenarios at $\alpha=0.25$ and  $\alpha=0.5$. %
We again use a solution time limit of 30 minutes.

Fig. \ref{fig:solve_times} shows the distribution of solution times for different OPS formulations, with a logarithmic y-axis. The results show that solution times for NF- and DC-OPS are generally similar and at least an order of magnitude lower than for SOC-OPS and several orders of magnitude lower than for the AC-OPS, validating the results from the IEEE 14 bus system. We also see that the solution times tend to increase as the systems size increases, making AC- and SOC-OPS too time consuming for the larger test cases. %
The difference in solve time when $\alpha \!\!=\!\! 0.25$ and $\alpha \!\!=\!\! 0.5$ becomes significant in scenarios with more than 24 buses, with $\alpha \!\! = \!\! 0.25$ leading to more time consuming problems.

\section{Conclusion} \label{sec:conclusion}

In this paper, we evaluate the trade-off in solution quality and solution time when modeling the OPS problem with the NF, DC, SOC and AC power flow formulations. To assess the solution quality for the NF, DC and SOC solutions, we find an AC-feasible power-flow solution given the de-energization decisions from the respective OPS problems. 

We find that solving the OPS problem with the linear NF and DC power flow formulations tend to significantly overestimate the amount of load that can be served in a given network configuration, thus resulting in solutions with high levels of load shed once the de-energization solutions are checked for AC feasibility. 
In comparison, the SOC-based OPS problems tend to more accurately assess the level of load shed, leading to better de-energization decisions. In our case study, SOC-OPS solutions even outperform the AC-OPS solutions, because the AC-OPS becomes stuck at local optima. 
In terms of solution time, NF and DC based formulations are orders of magnitude faster than SOC and AC power flow. 

Overall, results indicate that current solver technology forces us to pick between low solution quality with DC or NF power flow, or long solution time with SOC or AC power flow. Future work is needed to devise algorithms that can produce high quality solution in reasonable computational time.

\bibliographystyle{IEEEtran}
\bibliography{IEEEabrv,references}

\appendix
\edit{The variable bounds for variables $W^I_{ij}$ and $W^R_{ij}$ depend on the voltage angle difference bound for a line. In particular, the calculation changes if the minimum angle difference is greater than zero, or if the maximum angle difference is less than zero. Algorithm \ref{alg:bounds} shows the calculation of these variables bounds.}
\begin{algorithm}[h]
    \caption{Bounds calculation for $W^{R}_{ij}$ and $W^{I}_{ij}$}
        \begin{algorithmic}[1]\onehalfspacing
        \IF{$\boldsymbol{\underline{\theta_{ij}}} \ge 0$}
            \STATE $\boldsymbol{\overline{W}}^R_{ij}  = \boldsymbol{\overline{V}}_i  \boldsymbol{\overline{V}}_j  cos(\boldsymbol{\underline{\theta_{ij}}})$
            \STATE $\boldsymbol{\underline{W}}^R_{ij} = \boldsymbol{\underline{V}}_i \boldsymbol{\underline{V}}_j cos(\boldsymbol{\overline{\theta_{ij}}})$
            \STATE $\boldsymbol{\overline{W}}^I_{ij}  = \boldsymbol{\overline{V}}_i  \boldsymbol{\overline{V}}_j  sin(\boldsymbol{\overline{\theta_{ij}}})$
            \STATE $\boldsymbol{\underline{W}}^I_{ij} = \boldsymbol{\underline{V}}_i \boldsymbol{\underline{V}}_j sin(\boldsymbol{\underline{\theta_{ij}}})$
        \ELSIF{$\boldsymbol{\overline{\theta_{ij}}} \le 0$}
            \STATE $\boldsymbol{\overline{W}}^R_{ij}  = \boldsymbol{\overline{V}}_i  \boldsymbol{\overline{V}}_j  cos(\boldsymbol{\overline{\theta_{ij}}})$
            \STATE $\boldsymbol{\underline{W}}^R_{ij} = \boldsymbol{\underline{V}}_i \boldsymbol{\underline{V}}_j cos(\boldsymbol{\underline{\theta_{ij}}})$
            \STATE $\boldsymbol{\overline{W}}^I_{ij}  = \boldsymbol{\underline{V}}_i \boldsymbol{\underline{V}}_j sin(\boldsymbol{\overline{\theta_{ij}}})$
            \STATE $\boldsymbol{\underline{W}}^I_{ij} = \boldsymbol{\overline{V}}_i  \boldsymbol{\overline{V}}_j  sin(\boldsymbol{\underline{\theta_{ij}}})$
        \ELSE
            \STATE $\boldsymbol{\overline{W}}^R_{ij}  = \boldsymbol{\overline{V}}_i  \boldsymbol{\overline{V}}_j  $
            \STATE $\boldsymbol{\underline{W}}^R_{ij} = \boldsymbol{\underline{V}}_i \boldsymbol{\underline{V}}_j \min \{cos(\boldsymbol{\underline{\theta_{ij}}}), cos(\boldsymbol{\overline{\theta_{ij}}})\}$
            \STATE $\boldsymbol{\overline{W}}^I_{ij}  = \boldsymbol{\overline{V}}_i  \boldsymbol{\overline{V}}_j  sin(\boldsymbol{\overline{\theta_{ij}}})$
            \STATE $\boldsymbol{\underline{W}}^I_{ij} = \boldsymbol{\overline{V}}_i  \boldsymbol{\overline{V}}_j  sin(\boldsymbol{\underline{\theta_{ij}}})$
        \ENDIF
    \end{algorithmic}  \label{alg:bounds}
\end{algorithm}

\end{document}